# A variable service rate queue model for hub median problem

Zaniar Ardalan, Sajad Karimi

*Systems and Industrial Engineering department, University of Arizona, Tucson, Arizona, USA*
*Department of System Science and Industrial Engineering, Binghamton University, Binghamton, New York, USA*

**Abstract**

Hub location problems have multiple applications in logistic systems, airways industry, supply chain network design, and telecommunication. In the hub location problem, a number of nodes should be selected as the hub nodes to act as the main distributors and other nodes are connected together by these hubs. The input flow to the hub nodes is very large so more often we face the congestion in the hub nodes that causes disturbances in the whole system. Also, we have different service rates which is another cause of disturbance and should be addressed by the models. This paper addresses these issues by providing a model that prevents congestion in the system. We incorporated queuing system in the p-median hub location problem by considering multiple server options and different service rates. CAB dataset (contains 25 US cities) was used in the implementation and our findings show the big impact of considering congestion on the hub location network design.

*Keywords: hub location problem, p-hub median, queuing system, congestion, variable service rate,*

## 1. Introduction

Hub location problem is to design a network consisting of multiple locations that have inflow and outflow to each other. The goal is to design a network that minimizes the cost of establishing routes between locations, cost of distribution and delivery. Therefore, some locations will be chosen as hub locations and other locations will be connected to these locations.

In other words, mutual transportation flow between nodes should be established. A number of nodes selected as hub among all. These nodes act as the interface. Flow between any two nodes (hub nodes or non- hub nodes) must be connected through the hub nodes not using intermediate node. The number of arcs that we must pass is at least one arc and utmost two arcs.

In the literature four major types of hub location problems exist. 1. capacitated and uncapacitated hub location problem, 2. p-hub median problem, 3. p-hub center problem, 4. and hub covering location problems. The capacity of hub nodes may be capacitated/limited (LHLP) or uncapacitated/unlimited (UHLP). p-hub median problem (pHMP tries to minimize the total cost of transportation in the network by locating p hubs in an optimal manner. In this case, the number of hubs is given. Minimization of the farthest distance by finding the optimal location of p hubs and the allocation of non-hub nodes to the hubs is the main goal in the p-hub center problem (pHCP). Number of hubs is unknown in the hub covering location problem (HCLP). Locating hubs and allocating non-hub nodes to hub nodes are main problems here. It also contains cover constraints, which means for each hub we have limitations on the number of non-hub nodes it can serve.Hub problem in general can be categorized into two different types, single and multiple allocations problems. In the former, each non-hub node should be allocated to exactly one hub. In the latter type, a non-hub node can be allocated to more than one hub.

## 2. Literature Review

three coverage criteria was proposed Campbell [1]. The origin-destination pair *(i, j)* is covered by hubs *k* and *m* if:

1. The cost from *i* to *j* via *k* and *m* does not surpass a given quantity.
2. The cost of each arc from i to j via k and m does not surpass a given quantity.
3. Each of the origin–hub and hub–destination arcs satisfies different given values.

The multiple allocation and the uncapacitated single allocation problems were studied decades ago and Campbell [1] proposed the first mathematical model for the multiple allocation problem. Later, O'Kelly [2], developed multiple models for hub location problems.

They modeled the organization of a single and two hub networks. Klincewicz [3] proposed an effective algorithm for the uncapacitated hub location problem. Skorin-Kapov et al. [4], Ernst and Krishnamoorthy [5], and Mayer and Wagner [6] proposed multiple advances to the hub location problem, while Hamacher et al. [7] and Canovas et al. [8]studied the problem in the polyhedral manner. More recent advances have been presented by Marin [9] and Canovas et al. [10]. O'Kelly [11], Klincewicz [12], Skorin-Kapov et al. [4], Aykin [13] and Ernst and Krishnamoorthy [5] studied single allocation and Aykin [14], Ebery et al. [15], Campbell, [1] and [Boland et al. [16], and by Marin [17] studied the capacitated multiple allocation problem. The capacitated single allocation problem has also been studied by Ernst and Krishnamoorthy [19], Labbé et al. [18],Contreras et al. [19], and Contreras et al. [20], and by Contreras [21]. Campbell et al. [22] and Alumur and Kara [23] will be extensive sources for eager readers.

Riedi et al. [26] studied Multi-scale queuing (MSQ) and presented a model for that.. Ashour and Le-Ngoc [28], developed an MSQ model for variable-service rate multi-scale queuing (VS-MSQ) in order to evaluate priority queues. They presented an analytical framework to estimate the length of the queue and delay survivor functions for a priority queuing system with varying service rate. Considering congestion in the network was studied by Marianov and Serra [29]. They proposed a mathematical model for the network in order to find optimal locations of hub nodes. The most congested airports were modeled as M/D/c queuing systems. . They linearized the probabilistic constraint and solved the model using tabu search. Elhedhli and Xiaolong [30] model the congestion effect at a certain hub using a convex cost function that increases exponentially as more flows are directed through that hub. Mohammadi et al. [31] modeled Hubs, which are the most crowded parts of the network, as M/M/c queuing systems. In the literature of the hub location problem, multi-server hubs with different service rates are not addressed and this was our main motivation to incorporate this real world characteristic of the problem into the model. Rahimi et al. [33] incorporated congestion with uncertainty in the hub network design where they had multiple objective functions. Zhalechian et al. [34] incorporated social responsibility and congestion in the hub and spoke problem. Khodemoni-yazdi et al. [35] studied hierarchical hub location problems where they had two objective functions. They incorporated the queuing system in the paper. Sadeghi et al. [36]modeled p-hub covering problems considering budget for travel time. Zhalechian et al. [37] presented a new mathematical model to incorporate noise pollution of the transportation facilities.

In this paper we obtain the distribution rate of each node. The service capacity of each node should be studied to determine the probability distribution service rates. Given these rates and the existing relations in the theory of queues, queue length, waiting time in queue, time in the system are achieved. Because each node has the role of collector and distributor so each hub has a service system that includes one or more servers. Given the arrival rate and service rate of each hub, the system can be analyzed as a queuing system. In the posed system, each hub node is composed of several servers. The inflow to the hub node will be handled by multiple servers. Each of these servers has their own service rate. Service rate of servers depends on the system status and how many customers are present in the system.

For example, consider a case that cities are our nodes and they should be connected to each because of their mutual demands and supply. In this case hub cities can have multiple servers. In this case we want to use the hub median model with multiple allocation modes (Figure 1). If a city is selected as one of the hubs, we have to pay the cost of establishing some facilities. We also have transportation costs. The number of servers in each city is already known. We send $\beta_l$ percent of loads to $l^{th}$ server. In the cargo service centers, service rate varies exponentially. When congestion arises, managers can make decisions in order to speed up the process of service to facilitate the process. We want the number of customers within the queue to be less than a given value. The output of each server is identical in every respect. Quality service is similar in each of the servers. As a result, the output is similar for all hubs. And the outputs are ready to be sent to their final destination.

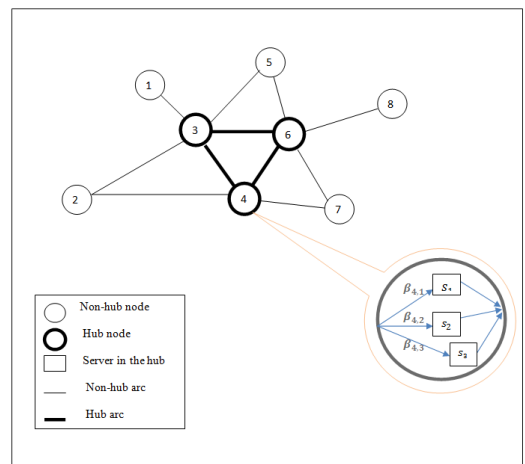

Fig.1. The schematic form of the proposed approach

## 2. Modeling:

This paper addresses a hub location problem based on the p-hub median model that considers the limited number of customers in the queue for each hub and the objective is minimizing the total transportation cost. P-hub median models use three-part paths for collection, transfer and distribution. Campbell [1] considers the cost for an origin–destination path as follows:

$C_{ijkm} = d_{ik} + \alpha d_{km} + d_{mj}$

Where $d_{km}$ is the distance between two hubs (k and m), and $d_{ik}$ and $d_{mj}$ are the distance between hub and non-hub nodes. The parameter $\alpha$ ($0 \leq \alpha \leq 1$) represents the discount corresponding to hub arcs to reflect the lower transportation cost due to the higher transportation scales.

The variables and parameters of the mathematical model are as follows:

Decision variables:

$X_{ijkm}$ : the fraction of flow from node $i$ (origin) to node $j$ (destination) from path $i$-$k$-$m$-$j$

$Y_k : \begin{cases} 1 & \text{if node } i \text{ is allocated to hub } k; \\ 0 & \text{otherwise,} \end{cases}$

Parameters:

$C_{ijkm}$: the transportation cost from node $i$ to node $j$ from path $i$-$k$-$m$-$j$

$W_{ij}$ : the flow (e.g. volume of freight) to be transported from node $i$ to node $j$.

$F_k$ : the fixed cost of opening a hub at node $k$.

$F_{km}$: the fixed cost related to the establishment of an arc between hub node $k$ and $m$.

$P$: the number of hubs.

$\theta_{q,kl}$: the maximum acceptable value for the probability of an excessive queue length at $l^{th}$ server of hub $k$.

$b_k$ : maximum acceptable queue length at hub $k$.

$$Min \sum_{i=1}^{n-1} \sum_{j=i+1}^{n} \sum_{k=1}^{n} \sum_{m=1}^{n} C_{ijkm} W_{ij} X_{ijkm} + \sum_{k=1}^{n} F_k Y_k$$
$$+ \sum_{k=1}^{n-1} \sum_{m=k+1}^{n} \frac{1}{2} F_{km} X_{kmkm} \quad (1)$$

Subject to:

$$\sum_{k=1}^{n} \sum_{m=1}^{n} X_{ijkm} = 1 \quad \forall (i,j,k,m) \in xfeas, i < j, \quad (2)$$

$$\sum_{k=1}^{n} Y_k = P \quad (3)$$

$$X_{ijkk} + \sum_{m \neq k} (X_{ijkm} + X_{ijmk}) \leq Y_k \quad \forall (i,j,k,m) \in xfeas, i < j, \quad (4)$$

$P$ [queue length at $l^{th}$ server of hub $k \geq b_{kl}$] $\leq \theta_{q,kl}$
$\forall (k, l)$ (5)

$X_{ijkm} \geq 0 \quad \forall (i,j,k,m) \quad (6)$

$Y_k \in \{0,1\} \quad \forall (k) \quad (7)$

The objective function (1) computes the total cost including the variable costs (collection, distribution and transportation) and fixed costs (hub and arc establishment). Constraint (2) ensures that each origin–destination flow goes through one or two hub nodes. Constraint (3) guarantees that exactly $p$ hubs are selected. Constraint (4) ensures that hubs are accessible for all of the transporters. Constraint (5) forces that the probability of exceeding the maximum acceptable trucks waiting in a queue in $l^{th}$ server of hub $k$ is less than or equal to $\theta_{q,kl}$. Constraints (6) and (7) define the sign of variables.

It is assumed that the service rate is in each server is a function of system status. Therefore, by changing the system status, the service rate is changed. In this model we assume that the relationship between service rate system status is as follows:

$\mu_n = n^c . \mu \quad (8)$

Where

$\mu$ is the average system rate when one truck is in the system.

$\mu_n$ is the average system rate when $n$ trucks are in the system.

Also, it is assumed that the input rate is according to Poisson process with parameter $\lambda$. Thus, the arrival rate is as follow:

$$\lambda_n = \lambda \qquad (9)$$

As mentioned before, each node consists of some servers and the servers are independent. The customers are allowed to choose each server on each node and wait in its queue. The arrival rate of each node is $\lambda_k$, therefore, according to the Poison's rule, the arrival rate of each server $\lambda_{kl}$ is calculated as follows:

$$\lambda_{kl} = \beta_{kl} \cdot \lambda_k \qquad \forall (k,l) \qquad (10)$$

And

$$\sum_{l=1}^{L} \beta_{kl} = 1 \qquad \forall (k) \qquad (11)$$

We need to rewrite the probabilistic equation (5) in an analytic form to make the mathematical model solvable. In order to obtain a deterministic linear constraint corresponding to this equation, we define the $p_s$ as steady-state probability of existing $s$ customers in the system with one server. Then, Eq. (5) can be written as:

$$\sum_{s=b_{kl}+2}^{\infty} p_s \leq \theta_{q,kl} \quad or \quad 1 - \sum_{s=0}^{b_{kl}+1} p_s \leq \theta_{q,kl} \qquad (12)$$

The left-hand side of the first form is related to the probability the probability of exceeding the maximum acceptable trucks waiting in a queue in $l^{th}$ server of hub $k$.

The second form has the same concept of the first form but it subtracts the probability of being less than upper bund trucks in queue from one.

In order to derive the expression equivalent to $p_s$, we do the following steps:

Assuming that the arrival rate and service rate of the system are $\lambda$ and $\mu$ respectively, the arrival and service rate of each state are:

$$\lambda_n = \lambda$$

$$\mu_n = n^c \cdot \mu$$

Also, in this model (an exponential model with variable service rates), the probability of being no truck in the hubs is:

$$p_0 = \left[1 + \sum_{n=1}^{\infty} \frac{\left(\frac{\lambda}{\mu}\right)^n}{(n!)^c}\right]^{-1} \qquad (13)$$

The probability of existing $n$ customer in the system is:

$$p_n = \left(\frac{\lambda}{\mu}\right)^n \frac{1}{(n!)^c} p_0 \qquad (14)$$

As there are infinite numbers of terms in eq (13), we propose a method to simplify computation of $p_0$. In this method, we define $M$ as a big number and replace the $\infty$ to $M$. Thus, the $p_0$ can be written as:

$$p_0 = \left[1 + \sum_{n=1}^{M} \frac{\left(\frac{\lambda}{\mu}\right)^n}{(n!)^c}\right]^{-1} \qquad (15)$$

The more quantity of $M$, the more accuracy in computing $p_0$. We define $\varepsilon$ as a minimum required accuracy for $p_0$.

$$\left[1 + \sum_{n=1}^{s-1} \frac{\left(\frac{\lambda}{\mu}\right)^n}{(n!)^c}\right]^{-1} - \left[1 + \sum_{n=1}^{s} \frac{\left(\frac{\lambda}{\mu}\right)^n}{(n!)^c}\right]^{-1} < \varepsilon \qquad (16)$$

Then we consider $M=s$ to compute $p_0$.

$$p_0 = \left[1 + \sum_{n=1}^{M} \frac{\left(\frac{\lambda}{\mu}\right)^n}{(n!)^c}\right]^{-1} \qquad (17)$$

Thus, the constraint (5) is:

$$\sum_{n=0}^{b_{kl}+1} \left(\frac{\lambda}{\mu}\right)^n \frac{1}{(n!)^c} \left[1 + \sum_{n=1}^{M} \frac{\left(\frac{\lambda}{\mu}\right)^n}{(n!)^c}\right]^{-1} \geq 1 - \theta_{q,kl} \qquad (18)$$

In this case, the location of hubs which is correspond to the variables $y_k$ will be determined after solving the model. On the other hand, the arrival rate to the hub $k$ is as follows:

$$\lambda_k = \left(\sum_{i=1}^{n}\sum_{j=1}^{n}\sum_{m=1}^{n} a_{ij} x_{ijkm} + \sum_{i=1}^{n}\sum_{j=1}^{n}\sum_{m=1}^{n} a_{ij} x_{ijmk}\right) (19)$$

and the arrival rate to the $l^{th}$ server of hub $k$ is given by:

$$\lambda_{kl} = \beta_{kl} \times \left(\sum_{i=1}^{n}\sum_{j=1}^{n}\sum_{m=1}^{n} a_{ij} x_{ijkm} \right.$$
$$\left. + \sum_{i=1}^{n}\sum_{j=1}^{n}\sum_{m=1}^{n} a_{ij} x_{ijmk}\right) \qquad (20)$$

Where $a_{ij}$ is the peak hour transportation flow from node $i$ to node $j$ by using hub $k$. Note that the total arrival rate to the hub $k$ is the sum of arrival rates from origin $i$ and hub $m$.

According to Marianov and Serra [20], we can numerically solve Eq. (18) and find the maximum value of variable $\lambda$ which is indicated by $\lambda_{max}$. All the possible values for the $\lambda$ that are less than $\lambda_{max}$ can satisfy the Eq. (18). Therefore, the Eq. (18) is equivalent to the following equation:

$$\lambda \leq \lambda_{max}$$

If the value of $\lambda_{max}$ is computed for each server of node $k$ (assuming that there is a difference between nodes in terms of service time), we could write $\lambda \leq \lambda_{max,kl}$ as follows:

$$\beta_{kl} \times \left( \sum_{i=1}^{n} \sum_{j=1}^{n} \sum_{m=1}^{n} a_{ij} x_{ijkm} + \sum_{i=1}^{n} \sum_{j=1}^{n} \sum_{m=1}^{n} a_{ij} x_{ijmk} \right) \leq \lambda_{max,kl} \qquad (21)$$

To solve the model, the following steps (that are proposed by Marianov and Serra [20]) are followed:

1. Estimating the service rate for each server of hubs.
2. Finding the $\lambda_{max,kl}$ for each server of hubs.
3. Solving the following model by using the value of $\lambda_{max,kl}$.

$$Min \sum_{i=1}^{n-1} \sum_{j=i+1}^{n} \sum_{k=1}^{n} \sum_{m=1}^{n} C_{ijkm} W_{ij} X_{ijkm} + \sum_{k=1}^{n} F_k Y_k + \sum_{k=1}^{n-1} \sum_{m=k+1}^{n} \frac{1}{2} F_{km} X_{kmkm} \qquad (1)$$

Subject to:

$$\sum_{k=1}^{n} \sum_{m=1}^{n} X_{ijkm} = 1 \quad \forall (i,j,k,m) \in xfeas, i < j, \qquad (2)$$

$$\sum_{k=1}^{n} Y_k = P \qquad (3)$$

$$X_{ijkk} + \sum_{m \neq k} (X_{ijkm} + X_{ijmk}) \leq Y_k \quad \forall (i,j,k,m) \in xfeas, i < j, \qquad (4)$$

$$\beta_{kl} \times \left( \sum_{i=1}^{n} \sum_{j=1}^{n} \sum_{m=1}^{n} a_{ij} x_{ijkm} + \sum_{i=1}^{n} \sum_{j=1}^{n} \sum_{m=1}^{n} a_{ij} x_{ijmk} \right) \leq \lambda_{max,kl} \qquad (21)$$

$$X_{ijkm} \geq 0 \quad \forall (i,j,k,m) \qquad (6)$$

$$Y_k \in \{0,1\} \quad \forall (k) \qquad (7)$$

### 3. Numerical example:

To evaluate our proposed approach, we conduct a numerical experiment on the CAB dataset that includes the distance and flow (intercity passengers) among 25 cities in the United States [32]. Figure 2 show these cities. Several instances are designed based on different value of saving transportation rates ($\alpha$=0.2, 0.5 and 0.8), number of hubs (p=4, 8 and 12) and allowable number of customers in queue (b= 5, 20).

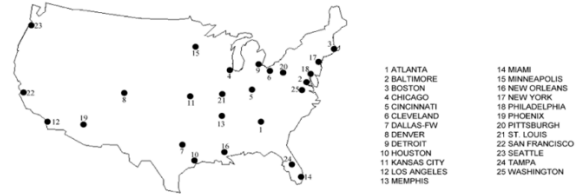

Fig.2. CAB data set

In these problems, a fixed cost is considered for establishment of hub nodes and hub-to-hub arcs are considered. We assume that each hub includes three servers with different service rates. The service rates of each server is a function of the number of people in the queue ($c = 0.2$, $\mu_{kln} = \mu_{kl} * n^{0.2}$). It is assumed that the customers enter to the first, second and third server with $\beta_{k1}$ =0.5, $\beta_{k2}$=0.2 and $\beta_{k3}$=0.3 respectively. For each server, the value of $\lambda_{kl}$ is calculated based on the $\theta_{q,kl} = 0.95$, and imported to the model. We code the model in GAMS 22.2 and solve on a PC with CI3, CPU 2.10 GHz, RAM 4 GB. The solution results (the value of objective function and the nodes selected as hubs) corresponding to each problem are represented in Table 1. According to the table 1, by increasing the value of $\alpha$, the total cost increases. Also, it is observed that the total cost for p=8 is less than the total cost for p=4 or p=12. This means that by increasing the number of hubs from 4 to 8, the transportation cost is reduced because the distances between non-hubs and hubs are reduced and the number of arcs with discounted fee is increased. In this situation, the transportation cost saving is more than the establishment costs for extra hubs. Conversely, by increasing the number of hubs from 8 to 12, the establishment cost would be higher than the transportation cost saving. Therefore, in these problems, the network with 8 hubs is more efficient than network with 4 or 12 hubs in terms of the total cost. Important point that is derived from the solution is that if constraint of allowable queue length altered from 5 to 20, total cost decreases because the problem is less restricted. In more detail, in the cases with lower allowable queue length, the model is forced to reduce the number of people in the queue by applying the hubs with higher service rates (that have higher establishment fixed cost).

| α | p | $b_{kl}$ | objective | Hub nodes |
|---|---|---|---|---|
| 0.2 | 4 | 5 | 5332 | 8, 14, 15, 17 |
| | | 20 | 4335 | 1, 8, 17, 20 |
| | 8 | 5 | 3361 | 2, 4, 6, 12, 13, 14, 17, 23 |
| | | 20 | 3061 | 2, 4, 6, 12, 13, 14, 17, 23 |
| | 12 | 5 | 3412 | 1, 2, 4, 6, 8, 13, 14, 15, 17, 21, 22, 23 |
| | | 20 | 3344 | 1, 2, 4, 6, 7, 13, 14, 16, 17, 19, 21, 23 |
| 0.5 | 4 | 5 | 5730 | 8, 14, 15, 17 |
| | | 20 | 4825 | 4, 8, 17, 20 |
| | 8 | 5 | 4149 | 1, 4, 12, 13, 14, 17, 20, 23 |
| | | 20 | 4037 | 1, 4, 12, 13, 14, 17, 20, 23 |
| | 12 | 5 | 4269 | 1, 2, 4, 6, 13, 14, 16, 17, 19, 21, 22, 23 |
| | | 20 | 4252 | 1, 2, 4, 6, 13, 14, 16, 17, 19, 21, 22, 23 |
| 0.8 | 4 | 5 | 5831 | 8, 14, 15, 17 |
| | | 20 | 5105 | 4, 8, 17, 20 |
| | 8 | 5 | 4833 | 1, 4, 14, 16, 17, 19, 20, 23 |
| | | 20 | 4780 | 1, 4, 14, 17, 19, 20, 21, 23 |
| | 12 | 5 | 5084 | 1, 2, 4, 6, 13, 14, 15, 16, 17, 19, 21, 23 |
| | | 20 | 5080 | 1, 2, 4, 6, 13, 14, 15, 16, 17, 19, 21, 23 |

Table 1: The solution of designed problems

In some solutions it is observed that increasing allowable length of queue, does not change the selected hubs and just changes the assignment of non-hubs to hubs that leads to reduce the total cost (for example : p=8 , α=0.5 and $b_{kl}$ altered from 5 to 20).

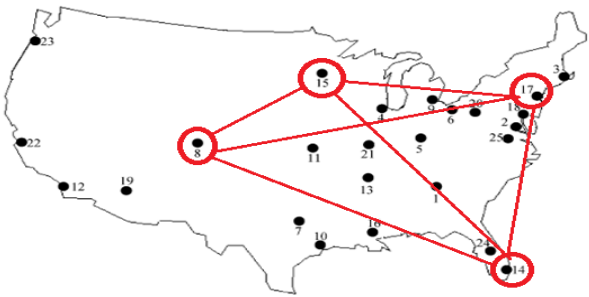

Fig.3 .hub nodes and hub arcs

The results of two computational solutions are shown in Figures 3 and 4. Figure 3 indicates the hub nodes and hub arcs for the problem with: p=4, α=0.2 and $b_{kl} = 5$ and Figure 4 indicates the hub nodes and hub arcs for the following parameters: p=4, α=0.2 and $b_{kl} = 20$.

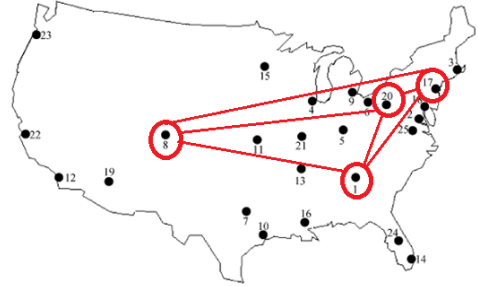

Fig.4 . hub nodes and hub arcs

## 5. Conclusion

In this paper, a p-hub median problem by considering the limitation for length of queue in the hubs was addressed. There are some servers in each hub that serve the customers (transportation facilities) to load, unload and distribute the freights. As the service rate of the servers is limited, a queue will be created in each server. The idea of this research is restricting the queue length in order to prevent the congestion in the hubs. The objective is minimizing the total transportation and fixed establishment costs regarding the aforementioned constraint. The model was formulated as a MILP model. A numerical experiment based on the CAB dataset was conducted to evaluate the proposed approach. The result of the numerical experiment indicated that the impact of restricting the queue length in the solutions. For example, decreasing the allowable queue length leads to increase the total cost due to the establishment of larger hubs with higher service rate and establishment cost. For the future works, developing a meta-heuristic for the proposed problem would be an interesting direction.